\documentclass[12pt]{amsart}
\usepackage{amsfonts}
\usepackage{amsmath,amssymb,amsthm}

\newcommand{\R}{\mathbb R}

\newtheorem{thm}{Theorem}[section]

\newtheorem{prop}[thm]{Proposition}
\theoremstyle{definition}
\newtheorem{defn}[thm]{Definition}
\theoremstyle{remark}

\newcommand{\ds}{\displaystyle}

\begin{document}

\title[INVARIANTS OF LINES ON SURFACES IN $\R^4$]
{INVARIANTS OF LINES ON SURFACES IN $\R^4$}

\author{Georgi Ganchev and Velichka Milousheva}
\address{Bulgarian Academy of Sciences, Institute of Mathematics and Informatics,
Acad. G. Bonchev Str. bl. 8, 1113 Sofia, Bulgaria}
\email{ganchev@math.bas.bg}
\address{Bulgarian Academy of Sciences, Institute of Mathematics and Informatics,
Acad. G. Bonchev Str. bl. 8, 1113, Sofia, Bulgaria; "L. Karavelov"
Civil Engineering Higher School, 175 Suhodolska Str., 1373 Sofia,
Bulgaria}
\email{vmil@math.bas.bg}



\begin{abstract}
Considering the tangent plane at a point to a surface in the four-dimensional
Euclidean space, we find an invariant of a pair of two tangents in this plane.
If this invariant is zero, the two tangents are said to be conjugate. When the
two tangents coincide with a given tangent, then we obtain the normal curvature
of this tangent. Asymptotic tangents (curves) are characterized by zero normal
curvature. Considering the invariant of the pair of a given tangent and its
orthogonal one, we introduce the geodesic torsion of this tangent. We obtain
that principal tangents (curves) are characterized by zero geodesic torsion.

The invariants $\varkappa$ and $k$ are introduced as the symmetric functions of
the two principal normal curvatures. The geometric meaning of the semi-sum
$\varkappa$ of the principal normal curvatures is equal (up to a sign) to the
curvature of the normal connection of the surface. The number of asymptotic
tangents at a point of the surface is determined by the sign of the invariant
$k$. In the case $k=0$ there exists a  one-parameter family of asymptotic lines,
which are principal. We find examples of such surfaces ($k=0$) in the class of
the general rotational surfaces (in the sense of Moore). The principal asymptotic
lines on these surfaces are helices in the four-dimensional Euclidean space.

\end{abstract}

\maketitle

\section{Preliminaries}

In \cite{GM1} we introduced a linear map $\gamma$ of Weingarten type
in the tangent space at any point of a surface in the four-dimensional
space $\R^4$. Analogously to the classical theory of surfaces in $\R^3$
the map $\gamma$ generates  second fundamental form $II$ on the surface.

In the present paper we introduce two invariants of a line on a surface
$M^2$ in $\R^4$ in terms of an invariant $\zeta_{\,g_1,g_2}$ of two
tangents $g_1$, $g_2$ at a point of $M^2$. Here we define
\emph{conjugate tangents} by the condition $\zeta_{\,g_1,g_2} = 0$.
We prove that the conjugacy in terms of $\zeta_{\,g_1,g_2}$ is the
conjugacy with respect to the second fundamental form. Thus we obtain
a geometric interpretation of the second fundamental form and
the Weingarten map of the surface.

Let $M^2: z = z(u,v), \, \, (u,v) \in {\mathcal D}$ (${\mathcal D}
\subset \R^2$) be a 2-dimensional surface in  $\R^4$ with tangent
space $T_pM^2 = {\rm span} \{z_u, z_v\}$ at any point $p \in M^2$,
and $E$, $F$, $G$ be the coefficients of the first fundamental form,
$W=\sqrt{EG-F^2}$. We choose an orthonormal normal frame field
$\{e_1, e_2\}$ of $M^2$ so that the quadruple $\{z_u, z_v, e_1, e_2\}$
is positive oriented in $\R^4$. Then the following derivative formulas
hold:
$$\begin{array}{l}
\vspace{2mm} \nabla'_{z_u}z_u=z_{uu} = \Gamma_{11}^1 \, z_u +
\Gamma_{11}^2 \, z_v
+ c_{11}^1\, e_1 + c_{11}^2\, e_2,\\
\vspace{2mm} \nabla'_{z_u}z_v=z_{uv} = \Gamma_{12}^1 \, z_u +
\Gamma_{12}^2 \, z_v
+ c_{12}^1\, e_1 + c_{12}^2\, e_2,\\
\vspace{2mm} \nabla'_{z_v}z_v=z_{vv} = \Gamma_{22}^1 \, z_u +
\Gamma_{22}^2 \, z_v
+ c_{22}^1\, e_1 + c_{22}^2\, e_2,\\
\end{array}$$
where $\Gamma_{ij}^k$ are the Christoffel's symbols and
$c_{ij}^k$, $i, j, k = 1,2$ are functions on $M^2$.

Denoting by $\sigma$ the second fundamental tensor of $M^2$, we have
$$\begin{array}{l}
\sigma(z_u,z_u)=c_{11}^1\, e_1 + c_{11}^2\, e_2,\\
[2mm]
\sigma(z_u,z_v)=c_{12}^1\, e_1 + c_{12}^2\, e_2,\\
[2mm] \sigma(z_v,z_v)=c_{22}^1\, e_1 + c_{22}^2\,e_2.
\end{array} \leqno{(1.1)}$$

The three pairs of normal vectors $\{\sigma(z_u,z_u),
\sigma(z_u,z_v)\}$, $\{\sigma(z_u,z_u), \sigma(z_v,z_v)\}$,
$\{\sigma(z_u,z_v), \sigma(z_v,z_v)\}$ form three parallelograms
with oriented areas
$$\Delta_1 = \left|%
\begin{array}{cc}
\vspace{2mm}
  c_{11}^1 & c_{12}^1 \\
  c_{11}^2 & c_{12}^2 \\
\end{array}%
\right|, \quad
\Delta_2 = \left|%
\begin{array}{cc}
\vspace{2mm}
  c_{11}^1 & c_{22}^1 \\
  c_{11}^2 & c_{22}^2 \\
\end{array}%
\right|, \quad
\Delta_3 = \left|%
\begin{array}{cc}
\vspace{2mm}
  c_{12}^1 & c_{22}^1 \\
  c_{12}^2 & c_{22}^2 \\
\end{array}%
\right|,$$ respectively. These oriented areas determine three
functions $\ds{L = \frac{2 \Delta_1}{W}, \,\, M = \frac{
\Delta_2}{W},\,\, N = \frac{2 \Delta_3}{W}}$, which change in the
same way as the coefficients $E, F, G$ under any change of the
parameters $(u,v)$.

Similarly to the theory of surfaces in $\R^3$, using the functions
$E$, $F$, $G$ and $L$, $M$, $N$ in \cite{GM1} we introduced the
linear map $\gamma$ in the tangent space at any point of $M^2$
$$\gamma: T_pM^2 \rightarrow T_pM^2,$$
defined by the equalities
$$\begin{array}{l}
\vspace{2mm}
\gamma(z_u)=\gamma_1^1z_u+\gamma_1^2z_v,\\
\vspace{2mm} \gamma(z_v)=\gamma_2^1z_u+\gamma_2^2z_v,
\end{array}$$
where
$$\displaystyle{\gamma_1^1=\frac{FM-GL}{EG-F^2}, \quad
\gamma_1^2 =\frac{FL-EM}{EG-F^2}}, \quad
\displaystyle{\gamma_2^1=\frac{FN-GM}{EG-F^2}, \quad
\gamma_2^2=\frac{FM-EN}{EG-F^2}}.$$

The linear map $\gamma$ of Weingarten type at the point $p \in
M^2$ is invariant with respect to changes of parameters on $M^2$
as well as to motions in $\R^4$. This implies that the functions
$$k = \frac{LN - M^2}{EG - F^2}, \qquad \varkappa =
\frac{EN+GL-2FM}{2(EG-F^2)}$$
are invariants of the surface $M^2$.

The invariant $\varkappa$ turns out to be the curvature of the normal
connection of the surface $M^2$ in $\R^4$.

As in the classical case, the invariants $k$ and $\varkappa$ divide
the points of $M^2$ into four types: flat, elliptic, parabolic and
hyperbolic. The surfaces consisting of flat points satisfy the conditions
$$k(u,v)=0, \quad \varkappa(u,v)=0, \qquad (u,v) \in \mathcal D,$$
or equivalently
$L(u,v)=0, \, M(u,v)=0, \, N(u,v)=0, \, (u,v) \in \mathcal D.$ These
surfaces are either planar surfaces (there exists a hyperplane
$\R^3 \subset \R^4$ containing $M^2$) or developable ruled surfaces
in $\R^4$ \cite{GM1}.

Further we consider surfaces free of flat points, i.e. $(L, M, N)
\neq (0, 0, 0)$.

\vskip 2mm Let $X = \lambda z_u + \mu z_v, \,\, (\lambda,\mu) \neq
(0,0)$ be a tangent vector at a point $p \in M^2$. The map $\gamma$
determines a second fundamental form of the surface $M^2$ at $p$
as follows:
$$II(\lambda,\mu) = - g(\gamma(X),X) =
L\lambda^2 + 2M\lambda\mu + N\mu^2, \quad \lambda,\mu \in \R.$$

As in the classical differential geometry of surfaces in $\R^3$ the second
fundamental form $II$ determines conjugate tangents at a point $p$ of $M^2$.
Two tangents $g_1: X_1 = \lambda_1 z_u + \mu_1 z_v$ and
$g_2: X_2 = \lambda_2 z_u + \mu_2 z_v$ are said to be
\textit{conjugate tangents}  if $II(\lambda_1, \mu_1; \lambda_2, \mu_2) = 0$, i.e.
$$L\lambda_1 \lambda_2 + M (\lambda_1 \mu_2 +\lambda_2 \mu_1) + N\mu_1 \mu_2 = 0.$$

In the next section we shall introduce conjugate tangents in a geometric way.

\section{Invariants of a tangent on a surface in $\R^4$}

Let $g$ be a tangent at the point $p \in M^2$ determined by the vector
$X =  \lambda z_u + \mu z_v$. We consider the  map
$\sigma_g: T_pM^2 \rightarrow (T_pM^2)^{\bot}$, defined by
$$\sigma_g(Y) = \ds{\sigma\left(\frac{\lambda z_u +
\mu z_v}{\sqrt{I(\lambda, \mu)}}, \, Y\right)}, \quad Y \in T_pM^2. \leqno{(2.1)}$$
Obviously $\sigma_g$ is a linear map, which does not depend on the choice of
the vector $X$ collinear with $g$. Using (1.1) and (2.1) we obtain the
following decomposition of the normal vectors $\sigma_g(z_u)$ and $\sigma_g(z_v)$:
$$\begin{array}{l}
\vspace{2mm}
\sigma_g(z_u)=\ds{\frac{\lambda\, c_{11}^1 + \mu \,c_{12}^1}
{\sqrt{I(\lambda, \mu)}}\, e_1 + \frac{\lambda\, c_{11}^2 +
\mu \,c_{12}^2}{\sqrt{I(\lambda, \mu)}}\, e_2},\\
\vspace{2mm}
\sigma_g(z_v)=\ds{\frac{\lambda\, c_{12}^1 + \mu \,c_{22}^1}
{\sqrt{I(\lambda, \mu)}}\, e_1 + \frac{\lambda\, c_{12}^2 +
\mu \,c_{22}^2}{\sqrt{I(\lambda, \mu)}}\, e_2}.
\end{array} \leqno{(2.2)}$$

Let $g_1: X_1 =  \lambda_1 z_u + \mu_1 z_v$ and $g_2: X_2 =
\lambda_2 z_u + \mu_2 z_v$ be two tangents at the point $p \in M^2$.
The oriented areas of the parallelograms determined by the pairs
of normal vectors $\sigma_{g_1} (z_u)$, $\sigma_{g_2} (z_v)$ and
$\sigma_{g_2} (z_u)$, $\sigma_{g_1} (z_v)$  are denoted by
$S(\sigma_{g_1} (z_u), \sigma_{g_2} (z_v))$, and
$S(\sigma_{g_2} (z_u), \sigma_{g_1} (z_v))$, respectively.
We assign the quantity $\zeta_{\,g_1,g_2}$ to the pair of tangents
$g_1$, $g_2$, defined by
$$\zeta_{\,g_1,g_2} = \ds{\frac{S(\sigma_{g_1} (z_u),\sigma_{g_2} (z_v))}{W} +
\frac{S(\sigma_{g_2} (z_u), \sigma_{g_1} (z_v))}{W}}. \leqno{(2.3)}$$

\begin{prop} The quantity $\zeta_{\,g_1,g_2}$ is invariant under any change of
the parameters on $M^2$.
\end{prop}
\noindent
\emph{Proof}:
Using equalities (2.2) we calculate that
$$\zeta_{\,g_1,g_2} = \ds{\frac{L \lambda_1 \lambda_2 +
M (\lambda_1 \mu_2 + \mu_1 \lambda_2) + N \mu_1 \mu_2}{\sqrt{I(\lambda_1, \mu_1)}
\sqrt{I(\lambda_2, \mu_2)}}} = \ds{\frac{II(\lambda_1, \mu_1; \lambda_2, \mu_2)}
{\sqrt{I(\lambda_1, \mu_1)} \sqrt{I(\lambda_2, \mu_2)}}}.$$

Now, let
$$\begin{array}{l}
\vspace{2mm}
u = u(\bar u,\bar v);\\
\vspace{2mm} v = v(\bar u,\bar v),
\end{array}
\quad (\bar u,\bar v) \in \bar{\mathcal D}, \,\, \bar{\mathcal D} \subset \R^2$$
be a smooth change of the parameters $(u,v)$
on $M^2$ with $J = u_{\bar u} \, v_{\bar v} - u_{\bar v} \,
v_{\bar u}\neq 0$. Then
$$\begin{array}{l}
\vspace{2mm}
z_{\bar u} = z_u \,u_{\bar u} + z_v \,v_{\bar u},\\
\vspace{2mm} z_{\bar v} = z_u \,u_{\bar v} + z_v \,v_{\bar v}
\end{array}$$
and $\bar E \bar G - \bar F^2=J^2\,(EG-F^2)$,
hence
$\bar W = \varepsilon J\,W, \,\, \varepsilon = {\rm sign} \, J$.
The functions $E$, $F$, $G$ and $L$, $M$, $N$ change as follows under the change
of the parametrization:
$$\begin{array}{ll}
\vspace{2mm}
\bar E=u_{\bar u}^2\,E+2\,u_{\bar u}v_{\bar u}\,F+v_{\bar u}^2\,G, & \quad \bar L=
\varepsilon (u_{\bar u}^2\,L+2\,u_{\bar u}v_{\bar u}\,M+v_{\bar u}^2\,N),\\
\vspace{2mm}
\bar F=u_{\bar u}u_{\bar v}\,E+(u_{\bar u}v_{\bar v}+v_{\bar u}u_{\bar v})\,F +
v_{\bar u}v_{\bar v}\,G, & \quad
\bar M = \varepsilon (u_{\bar u}u_{\bar v}\,L+(u_{\bar u}v_{\bar v}+
v_{\bar u}u_{\bar v})\,M + v_{\bar u}v_{\bar v}\,N)\\
\vspace{2mm}
\bar G=u_{\bar v}^2\,E+2\,u_{\bar v}v_{\bar v}\,F+v_{\bar v}^2\,G, & \quad
\bar N=\varepsilon (u_{\bar v}^2\,L+2\,u_{\bar v}v_{\bar v}\,M+v_{\bar v}^2\,N).
\end{array}$$
If $X =  \lambda z_u + \mu z_v = \bar \lambda z_{\bar u} + \bar \mu z_{\bar v}$, then
$\lambda = u_{\bar u} \bar \lambda + u_{\bar v} \bar \mu, \,\,\mu =
v_{\bar u} \bar \lambda + v_{\bar v} \bar \mu$.

The first and the second fundamental forms change as follows:
$$\bar I (\bar \lambda, \bar \mu) = I(\lambda, \mu); \qquad
\bar{II} (\bar \lambda, \bar \mu) = \varepsilon II(\lambda, \mu).$$
Hence,
$$\bar{\zeta}_{\,g_1, g_2} = \ds{\frac{\bar{II}(\bar \lambda_1, \bar \mu_1;
\bar \lambda_2, \bar \mu_2)}{\sqrt{I(\bar \lambda_1, \bar \mu_1)}
\sqrt{I(\bar \lambda_2,\bar \mu_2)}}}=
\ds{\frac{\varepsilon \,II(\lambda_1, \mu_1; \lambda_2, \mu_2)}
{\sqrt{I(\lambda_1, \mu_1)}
\sqrt{I(\lambda_2, \mu_2)}}} = \varepsilon \, \zeta_{\,g_1,g_2}.$$
Consequently, $\zeta_{\,g_1,g_2}$ is invariant (up to the orientation of the tangent space or the normal space).
\qed

\vskip 3mm
\begin{defn} \label{D:conjugate}
Two tangents $g_1: X_1 = \lambda_1 z_u + \mu_1 z_v$ and $g_2: X_2 = \lambda_2 z_u + \mu_2 z_v$ are said to
be \emph{conjugate tangents}   if $\zeta_{\,g_1,g_2} = 0$.
\end{defn}

Obviously, $\zeta_{\,g_1,g_2} = 0$ if and only if
$L\lambda_1 \lambda_2 + M (\lambda_1 \mu_2 +\lambda_2 \mu_1) + N\mu_1 \mu_2 = 0.$
Hence, the tangents $g_1$ and $g_2$ are conjugate according to Definition \ref{D:conjugate}
if and only if they are conjugate with respect to  the second fundamental form $II$.

\vskip 2mm
We shall assign two invariants $\nu_g$ and $\alpha_g$ to any tangent $g$ of the surface in the following way.
Let $g: X = \lambda z_u + \mu z_v$ be a tangent and $g^{\bot}$ be its orthogonal tangent, determined by the vector
$$X^{\bot} = \ds{-\frac{F \lambda + G \mu}{W} \, z_u + \frac{E \lambda + F \mu}{W} \, z_v}.  \leqno{(2.4)}$$

We define
$$\nu_g = \zeta_{\,g,g}; \qquad \alpha_g = \zeta_{\,g, g^{\bot}}. $$
We call $\nu_g$ the \emph{normal curvature} of the tangent $g$, and $\alpha_g$ - the \emph{geodesic torsion} of $g$.

The equality (2.3) implies that
$$\nu_g = \ds{2\frac{S(\sigma_g(z_u), \sigma_g(z_v))}{W}} = \ds{\frac{II(\lambda, \mu)}{I(\lambda, \mu)}}.$$
Hence, the normal curvature of the tangent $g$ is two times  the oriented area of the parallelogram determined
by the normal vectors $\sigma_g(z_u)$ and $\sigma_g(z_v)$. The invariant $\nu_g$ is expressed by the first and the second fundamental forms of the surface
in the same way as the normal curvature of a tangent in the theory of surfaces in $\R^3$.

Using (2.3) and (2.4) we get
$$\alpha_g = \ds{\frac{\lambda^2 (EM - FL) + \lambda \mu (EN - GL) + \mu^2(FN - GM)}{W I(\lambda, \mu)}}.$$
The last formula shows that $\alpha_g$ is expressed by the coefficients of the first and the second fundamental forms in the same way
as the geodesic torsion in the theory of surfaces in $\R^3$.

\vskip 3mm
\begin{defn} \label{D:asymptotic tangent}
A tangent $g: X = \lambda z_u + \mu z_v$ is said to be \emph{asymptotic} if it is self-conjugate.
\end{defn}

\begin{prop}
A tangent $g$ is asymptotic if and only if $\nu_g = 0$.
\end{prop}

\vskip 3mm
\begin{defn} \label{D:principal tangent}
A tangent $g: X = \lambda z_u + \mu z_v$ is said to be \emph{principal} if it is perpendicular to its conjugate.
\end{defn}

\begin{prop}
A tangent $g$ is principal if and only if $\alpha_g = 0$.
\end{prop}

\vskip 3mm
The equation for the asymptotic tangents at a point $p \in M^2$ is
$$L\lambda^2 + 2M \lambda \mu + N\mu^2 = 0.$$
If $p$ is an elliptic point of $M^2$ ($k>0$) then there are no asymptotic tangents through $p$;
if $p$ is a hyperbolic point ($k<0$)  then there are two asymptotic tangents passing through $p$, and if $p$ is a parabolic point ($k=0$)
then there is one asymptotic tangent through $p$.

A line $c: u=u(q), \; v=v(q); \; q\in \textrm{J} \subset \R$ on $M^2$ is
said to be an \textit{asymptotic line}  if its tangent at any point is asymptotic.

\vskip 2mm
The equation for the principal tangents at a point $p \in M^2$ is
$$\left|\begin{array}{cc}
E & F\\
[2mm] L & M \end{array}\right| \lambda^2+ \left|\begin{array}{cc}
E & G\\
[2mm] L & N \end{array}\right| \lambda \mu+
\left|\begin{array}{cc}
F & G\\
[2mm] M & N \end{array}\right| \mu^2=0.$$

A line $c: u=u(q), \; v=v(q); \; q\in \textrm{J} \subset \R$ on $M^2$ is
said to be a \textit{principal line} (a \textit{line of
curvature}) if its tangent at any point is principal. The surface
$M^2$ is parameterized by the principal lines if and
only if $F=0,\,\, M=0$.

The normal curvatures $\nu' = \ds{\frac{L}{E}}$ and $\nu'' =
\ds{\frac{N}{G}}$ of the principal tangents are said to be
\textit{principal normal curvatures} of $M^2$.
The invariants $k$ and $\varkappa$ of $M^2$ are expressed by the
principal normal curvatures $\nu'$ and $\nu''$ as follows:
$$k = \nu' \nu''; \qquad \varkappa = \frac{\nu' + \nu''}{2}.$$

\section{Examples of surfaces with $k = 0$} \label{S:Examples}

In this section we consider general (in the sense of C.
Moore) rotational surfaces in $\R^4$ whose meridians lie in  two-dimensional planes.
We shall find all such surfaces consisting of parabolic points.

Considering general rotations in $\R^4$, C. Moore introduced
general rotational surfaces \cite{M} \, (see also \cite{MW1,
MW2}). In the case when the meridians lie in two-dimensional planes
the general rotational surface can be parameterized as follows:
$$M^2: z(u,v) = \left( f(u) \cos\alpha v, f(u) \sin \alpha v, g(u) \cos \beta v, g(u) \sin \beta v \right);
\quad u \in \textrm{J} \subset \R, \,\,  v \in [0; 2\pi),$$
where $f(u)$ and $g(u)$ are smooth functions, satisfying $\alpha^2
f^2(u)+ \beta^2 g^2(u) > 0 , \,\, f'\,^2(u)+ g'\,^2(u) > 0,\, u
\in \textrm{J}$, and $\alpha, \beta$ are positive constants.

Each parametric curve $u = u_0 = const$ of $M^2$ is given by
$$c_v: z(v) = \left( a \cos \alpha v, a \sin \alpha v, b \cos \beta v, b \sin \beta v \right);
\quad a = f(u_0), \,\, b = g(u_0)$$ and its Frenet curvatures are
$\varkappa_{c_v} = \ds{\sqrt{\frac{a^2 \alpha^4 + b^2 \beta^4}{a^2 \alpha^2 + b^2 \beta^2}}}; \,\,
\tau_{c_v} = \ds{\frac{ab \alpha \beta (\alpha^2 - \beta^2)}
{\sqrt{a^2 \alpha^4 + b^2 \beta^4}\sqrt{a^2 \alpha^2 + b^2
\beta^2}}}; \,\, \sigma_{c_v} = \ds{\frac{\alpha \beta \sqrt{a^2
\alpha^2 + b^2 \beta^2}} {\sqrt{a^2 \alpha^4 + b^2 \beta^4}}}$.
Hence, when $\alpha \neq \beta$ each parametric curve $u =
const$ is a curve in $\R^4$ with constant curvatures (helix), and when
 $\alpha = \beta$ each parametric curve $u = const$ is a circle.

Each parametric curve $v = v_0 = const$ of $M^2$ is given by
$$c_u: z(u) = \left(\, A_1 f(u), A_2 f(u), B_1 g(u), B_2 g(u) \, \right),$$
where  $A_1 = \cos \alpha v_0, \, A_2 = \sin \alpha v_0, \,B_1 =
\cos \beta v_0, \,B_2 = \sin \beta v_0$. The Frenet curvatures of
$c_u$ are $\varkappa_{c_u} = \ds{\frac{|g' f'' - f' g''|}{(\sqrt{f'\,^2 + g'\,^2})^3}}; \,\, \tau_{c_u} = 0$.
Hence, $c_u$ is a plane curve with curvature $\varkappa_{c_u} =
\ds{\frac{|g' f'' - f' g''|} {(\sqrt{f'\,^2 + g'\,^2})^3}}$. So,
for each $v = const$ the parametric curves $c_u$ are congruent in
$\R^4$. These curves are called \textit{meridians} of $M^2$.

The tangent space of $M^2$ is spanned by the vector fields
$$\begin{array}{l}
\vspace{2mm}
z_u = \left(f' \cos \alpha v, f' \sin \alpha v, g' \cos \beta v, g' \sin \beta v \right);\\
\vspace{2mm} z_v = \left( - \alpha f \sin \alpha v, \alpha f \cos
\alpha v, - \beta g \sin \beta v, \beta g \cos \beta v \right).
\end{array}$$
Hence, the coefficients of the first fundamental form are $E =
f'\,^2(u)+ g'\,^2(u); \,\, F = 0; \,\, G = \alpha^2 f^2(u)+
\beta^2 g^2(u)$. We consider the following orthonormal tangent
frame field
$$\begin{array}{l}
\vspace{2mm}
x = \ds{\frac{1}{\sqrt{f'\,^2 + g'\,^2}}\left(f' \cos \alpha v, f' \sin \alpha v,
g' \cos \beta v, g' \sin \beta v \right)};\\
\vspace{2mm} y = \ds{\frac{1}{\sqrt{\alpha^2 f^2 + \beta^2
g^2}}\left( - \alpha f \sin \alpha v, \alpha f \cos \alpha v, -
\beta g \sin \beta v, \beta g \cos \beta v \right)};
\end{array}$$
and the following orthonormal normal frame field
$$\begin{array}{l}
\vspace{2mm}
n_1 = \ds{\frac{1}{\sqrt{f'\,^2 + g'\,^2}}\left(g' \cos \alpha v,
g' \sin \alpha v, - f' \cos \beta v, - f' \sin \beta v \right)};\\
\vspace{2mm} n_2 = \ds{\frac{1}{\sqrt{\alpha^2 f^2 + \beta^2
g^2}}\left( - \beta g \sin \alpha v, \beta g \cos \alpha v, \alpha
f \sin \beta v, - \alpha f \cos \beta v \right)}.
\end{array}$$
$\{x, y, n_1, n_2\}$ is a positive
oriented orthonormal frame field in $\R^4$.

We calculate the functions $c_{ij}^k, \,\, i,j,k = 1,2$ and find the coefficients $L$, $M$ and $N$ of the second
fundamental form of $M^2$:
$$L = \ds{\frac{2 \alpha \beta (g f' - f g') (g' f'' - f' g'')}{(\alpha^2 f^2
+ \beta^2 g^2) (f'\,^2 + g'\,^2)}}; \quad M = 0; \quad
N = \ds{\frac{- 2\alpha \beta (g f' - f g') (\beta^2 g f' -
\alpha^2 f g')}{(\alpha^2 f^2 + \beta^2 g^2) (f'\,^2 +
g'\,^2)}}.$$

Consequently, the invariants $k$, $\varkappa$ and the Gauss curvature $K$ of $M^2$
are expressed as follows:

$$k = \ds{\frac{- 4 \alpha^2 \beta^2 (g f' - f g')^2 (g' f'' - f' g'')
(\beta^2 g f' - \alpha^2 f g')}{(\alpha^2 f^2 + \beta^2 g^2)^3 (f'\,^2 + g'\,^2)^3}};$$

$$\varkappa =  \ds{\frac{\alpha \beta (g f' - f g')}
{(\alpha^2 f^2 + \beta^2 g^2)^2 (f'\,^2 + g'\,^2)^2} \,
\left( (\alpha^2 f^2 + \beta^2  g^2)(g' f'' - f' g'') -  (f'\,^2 +
g'\,^2) (\beta^2 g f' - \alpha^2 f g') \right)};$$

$$K =  \ds{\frac{(\alpha^2 f^2 + \beta^2  g^2)(\beta^2 g f'
- \alpha^2 f g')(g' f'' - f' g'') - \alpha^2 \beta^2 (f'\,^2 + g'\,^2) (g f' - f g')^2}
{(\alpha^2 f^2 + \beta^2 g^2)^2 (f'\,^2 + g'\,^2)^2} \,}.$$

Note that the surface $M^2$ is parameterized by the principal lines ($F = M = 0$).

\vskip 2mm Now we shall find the general rotational surfaces
with $k = 0$. Without loss of generality we assume that the
meridian $m$ is defined by $f = u; \,\, g = g(u)$. Then
$$k = \ds{\frac{ 4 \alpha^2 \beta^2 (g  - u g')^2 g'' (\beta^2 g  - \alpha^2 u g')}
{(\alpha^2 u^2 + \beta^2 g^2)^3 (1 + g'\,^2)^3}}.$$

The invariant $k$ is zero in the following three cases:

\vskip 2mm 1. $g(u) = a\,u,\; a = const \neq 0$. In this case $k =
\varkappa = K = 0$, and $M^2$ is a developable surface in $\R^4$.

\vskip 2mm 2. $g(u) = a\,u + b ,\; a= const \neq 0, b = const \neq
0$. In this case $k = 0$, but $\varkappa \neq 0$, $K \neq 0$.
Consequently, $M^2$ is a non-developable ruled surface in $\R^4$.

\vskip 2mm 3. $g(u) = \ds{ c\,u^{\frac{\beta^2}{\alpha^2}}}, \; c
= const \neq 0$. In case of $\alpha \neq \beta$ we get $k = 0$,
and the invariants $\varkappa$ and $K$ are given by
$$\varkappa = \ds{\frac{c^2 \alpha^{-5} \beta^3 (\beta^2 - \alpha^2 )^2 u^{2 \frac{\beta^2-\alpha^2}{\alpha^2}}}
{ \left(\alpha^2 u^2 + \beta^2 c^2 u^{2
\frac{\beta^2}{\alpha^2}}\right) \left(1 + c^2
\frac{\beta^4}{\alpha^4} u^{2\frac{\beta^2-\alpha^2}{\alpha^2}
}\right)^2}}; \quad
K = \ds{- \frac{c^2 \alpha^{-2} \beta^2 (\beta^2 - \alpha^2 )^2 u^{2 \frac{\beta^2}{\alpha^2}}}
{\left(\alpha^2 u^2 + \beta^2 c^2 u^{2
\frac{\beta^2}{\alpha^2}}\right)^2 \left(1 + c^2
\frac{\beta^4}{\alpha^4}
u^{2\frac{\beta^2-\alpha^2}{\alpha^2}}\right)}}.$$ Hence,
$\varkappa \neq 0$, $K \neq 0$. In this case the parametric lines
$u =const$ and $v = const$ are not straight lines. This is the most interesting
example of general rotational surfaces with $k = 0$.

Since $k=0$,  one of the principal normal curvatures of $M^2$ is zero ($\nu'' = \ds{\frac{N}{G} = 0}$).
Hence, the parametric $v$-lines of $M^2$ are asymptotic principal lines.
Moreover, these lines are helices in $\R^4$.

\vskip 2mm
\textbf{Acknowledgements:} The second author is partially supported by "L. Karavelov"
Civil Engineering Higher School, Sofia, Bulgaria under Contract No 10/2010.

\end{document}